\documentclass{article}

\usepackage{myarxiv}
\usepackage[utf8]{inputenc} 
\usepackage[T1]{fontenc}    
\usepackage{hyperref}       
\usepackage{url}            
\usepackage{booktabs}       
\usepackage{amsfonts}       
\usepackage{nicefrac}       
\usepackage{microtype}      
\usepackage{lipsum}

\usepackage{amssymb,amsmath,amsfonts,amsthm,enumerate}
\usepackage{subfigure}
\usepackage{graphicx}
\usepackage{epsfig}
\usepackage{float}
\usepackage[usenames]{color}
\usepackage{cite}
\usepackage{bm}
\usepackage{enumitem}
\setlist{nosep}

\newtheorem{lemma}{Lemma}[section]

\newtheorem{theorem}{Theorem}[section]

\def\BC{\mathbb C}

\def\BR{\mathbb R}

\def\cD{\mathcal D}

\def\cU{\mathcal U}

\def\rd{\mathrm d}

\def\supp{\mathrm{supp}}

\def\Ga{\Gamma}

\def\Om{\Omega}
\def\al{\alpha}
\def\be{\beta}
\def\ga{\gamma}
\def\de{\delta}

\def\ve{\varepsilon}

\def\ka{\kappa}
\def\la{\lambda}
\def\si{\sigma}
\def\vp{\varphi}
\def\om{\omega}
\def\f{\frac}

\def\ov{\overline}
\def\pa{\partial}

\def\tri{\triangle}

\title{Inverse Problems of Determining Sources of the Fractional Partial Differential Equations}

\author{
  Yikan Liu\\
  Graduate School of Mathematical Sciences\\
  The University of Tokyo\\
  3-8-1 Komaba, Meguro-ku, Tokyo 153-8914, Japan\\
  \texttt{ykliu@ms.u-tokyo.ac.jp}\\
  \And
  Zhiyuan Li\\
  School of Mathematics and Statistics\\
  Shandong University of Technology\\
  Zibo, Shandong 255049, China\\
  \texttt{zyli@sdut.edu.cn}\\
  \AND
  Masahiro Yamamoto\\
  Graduate School of Mathematical Sciences\\
  The University of Tokyo\\
  3-8-1 Komaba, Meguro-ku, Tokyo 153-8914, Japan\\
  \texttt{myama@ms.u-tokyo.ac.jp}\\
}

\begin{document}
\maketitle

\begin{abstract}
In this chapter, we mainly review theoretical results on inverse source problems for diffusion equations with the Caputo time-fractional derivatives of order $\al\in(0,1)$. Our survey covers the following types of inverse problems:
\begin{itemize}
\item determination of time-dependent functions in interior source terms
\item determination of space-dependent functions in interior source terms
\item determination of time-dependent functions appearing in boundary conditions
\end{itemize}
\end{abstract}

\MRsubject{35R11\and 26A33\and 35R30\and 65M32}

\section{Introduction}\label{sec-intro}

Let $\Om\subset\BR^d$ be a bounded domain with smooth boundary $\pa\Om$, $T>0$, $\al\in(0,1)$ and $x=(x_1,\ldots,x_d)\in\BR^d$, $\pa_t=\f\pa{\pa t}$, $\pa_j=\f\pa{\pa x_j}$. By $\tri:=\sum_{j=1}^d\pa_j^2$, we denote the usual Laplacian in space, and $\pa_t^\al$, $0<\al<1$ stands for the Caputo derivative in time:
\[
\pa_t^\al f(t):=\f1{\Ga(1-\al)}\int_0^t (t-s)^{-\al}\f{\rd f}{\rd s}(s)\,\rd s,
\]
where $\Ga$ denotes the gamma function.

Then we consider an initial-boundary value problem
\begin{equation}\label{eq-ibvp-u}
\begin{cases}
(\pa_t^\al-\tri)u(x,t)=g(x)\rho(t), & x\in\Om,\ 0<t<T,\\
u(x,0)=0, & x\in\Om,\\
u(x,t)=0, & x\in\pa\Om,\ 0<t<T.
\end{cases}
\end{equation}
The initial-boundary value problem governs the time evolution of density $u(x,t)$ at location $x$ and time $t$ of some substance such as contaminants. Here, $g(x)\rho(t)$ is an interior source term producing the substance in $\Om$, and the source term can be often assumed to be modeled in the form of the separation of variables. Here, $g$ and $\rho$ describe the spatial distribution of the source and the time evolution pattern, respectively. In the typical inverse problems for \eqref{eq-ibvp-u}, we are requested to determine $\rho(t)$ and/or $g(x)$ by extra data of solution $u(x,t)$ to \eqref{eq-ibvp-u}.

Mainly, we survey the following two types of inverse source problems.\smallskip

{\bf Inverse $t$-source problem}

Let $u$ satisfy \eqref{eq-ibvp-u}. We fix $x_0\in\Om$ arbitrarily. Provided that the spatial component $g$
in the source term is known, determine the temporal component $\rho$ in $(0,T)$ by the single point observation of $u$ in $\{x_0\}\times(0,T)$.\smallskip

{\bf Inverse $x$-source problem}

Let $u$ satisfy \eqref{eq-ibvp-u}.  We suppose that the temporal component $\rho$ in the source term is known. Determine the spatial component $g$ by
\begin{enumerate}
\item[(a)] the final time observation of $u$ in $\Om\times\{T\}$, or
\item[(b)] the partial interior observation of $u$ in $\om\times(0,T)$, where $\om\subset\Om$ is an arbitrarily chosen nonempty subdomain.
\end{enumerate}\smallskip
The above inverse $t$-source problem and the inverse $x$-source problem (a) with the final time observation have been well studied and many theoretical researches have been published for classical partial differential equations. As monograph, we should refer to Prilepko, Orlovsky, and Vasin \cite{POV00}. Moreover see for example Choulli and Yamamoto \cite{CY96,CY97,CY99}, Isakov \cite{I91}, Tikhonov and Eidelman \cite{TE94} and the references therein for related inverse problems for usual partial differential equations.

As for the fractional differential equations, we know that we can construct theories parallel to \cite{POV00}, and the works are now made continuously. This chapter is nothing but an incomplete survey for these works in full progress.

The identification of $\rho(t)$ fits, for example, in the cases of disasters of nuclear power plants, in which the source location can be assumed to be known but the decay of the radiative strength in time is unknown and important to be estimated. On the other hand, one example of the identification of $g(x)$ can be illustrated by the detection of illegal discharge of sewage, which is a serious problem in some countries. Also by such practical demands, the inverse source problems have been strongly required to be studied both theoretically and numerically.

Some papers surveyed later discuss the case where $-\tri$ is replaced by a general linear elliptic operator:
\[
-\sum_{i,j=1}^d\pa_j\left(a_{ij}(x,t)\pa_i u\right)+\sum_{j=1}^d b_j(x,t)\pa_j u+c(x,t)u.
\]
Moreover, we can similarly discuss the following generalized cases:
\begin{itemize}
\item Multi-term $\sum_{j=1}^m q_j\pa_t^{\al_j}u$ or distributed order $\int_0^t\mu(\al)\pa_t^\al u\,\rd\al$ in time (see, e.g., the chapter on ``Inverse problems of determining parameters of the fractional partial differential equations'' of this handbook). Here, $0<\al_m<\cdots<\al_1<1$, $q_1,\ldots,q_m$ are constants, and $\mu\in L^1(0,1),\ge0,\not\equiv0$.
\item$\rho(t)g(x,t)$ in the inverse $t$-source problem and $g(x)\rho(x,t)$ in the inverse $x$-source problem.
\item Other kinds of boundary conditions and inhomogeneous boundary values, and the whole space $\BR^d$.
\end{itemize}
However, in order to focus on the main topic, we choose a simple formulation \eqref{eq-ibvp-u}, which definitely captures
the essence of the inverse source problems.

We note that in the cases of $\rho(t)g(x,t)$ and $g(x)\rho(x,t)$, the inverse source problems are considered as linearized
problems of inverse coefficient problems (e.g., see the chapter on ``Inverse problems of determining coefficients of the fractional partial differential equations'' of this handbook). For example, we discuss the determination of $p(x)$ in
\[
(\pa_t^\al-\tri)y(x,t)=p(x)y(x,t),\quad x\in \Om,\ 0<t<T
\]
by extra data. Let $y$ and $z$ be the corresponding solutions respectively with the coefficients $p(x)$ and $q(x)$. Then, setting $u=y-z$, $g(x)=p(x)-q(x)$ and $\rho(x,t)=z(x,t)$, we reduce the inverse problem of determining a coefficient to an inverse source problem for
\[
(\pa_t^\al-\tri-p(x))u(x,t)=g(x)\rho(x,t),\quad x\in \Om,\ 0<t<T.
\]

Most of publications on inverse problems for fractional equations are concerned with the case of the Caputo derivatives in $t$ of orders $\al\in(0,1)$. Mainly, we discuss the Caputo derivative $\pa_t^\al$ with $0<\al<1$, although other kinds of fractional derivatives and/or $1<\al<2$ are also meaningful.

In this chapter, we mainly review theoretical results in the existing literature for inverse source problems with slight improvements time by time. We will at most provide key ideas and sketches of the proofs instead of detailed arguments.

The remainder of this chapter is organized as follows. In Section \ref{sec-pre}, we prepare the necessary ingredients for dealing with the inverse problems including the basic facts on forward problems of \eqref{eq-ibvp-u}. In Sections \ref{sec-ISP-t}--\ref{sec-ISP-x}, we review the inverse $t$- and $x$-source problems, respectively. In Section \ref{sec-re}, we survey related inverse problems of determining some functions in boundary conditions. Finally, we summarize the chapter with concluding remarks in Section \ref{sec-con}.

\section{Preliminaries}\label{sec-pre}

Henceforth, $L^2(\Om)$, $H^1_0(\Om)$, $H^2(\Om)$ are the usual $L^2$-space and the Sobolev spaces of real-valued functions
(e.g., Adams \cite{A75}): $L^2(\Om)=\{f;\ \int_\Om|f(x)|^2\,\rd x<\infty\}$, and let $(\,\cdot\,,\,\cdot\,)$ denote the scalar product: $(f,g)=\int_\Om f(x)g(x)\,\rd x$ for $f,g\in L^2(\Om)$. We define the Laplace operator $-\tri$ with the domain $\cD(-\tri)=H^2(\Om)\cap H^1_0(\Om)$. An eigensystem $\{(\la_n,\vp_n)\}_{n=1}^\infty$ of $-\tri$ is defined
by $\{\vp_n\}_{n=1}^\infty\subset\cD(-\tri)$ and $0<\la_1<\la_2\le\ldots$ such that $-\tri\vp_n=\la_n\vp_n$ and $\{\vp_n\}_{n=1}^\infty$ is a complete orthonormal system of $L^2(\Om)$.

For $\ga\ge0$, the fractional Laplacian $(-\tri)^\ga$ is defined as
\begin{align*}
\cD((-\tri)^\ga) & :=\left\{f\in L^2(\Om);\sum_{n=1}^\infty\left|\la_n^\ga(f,\vp_n)\right|^2<\infty\right\},\\
(-\tri)^\ga f & :=\sum_{n=1}^\infty\la_n^\ga(f,\vp_n)\vp_n.
\end{align*}
We know that $\cD((-\tri)^\ga)$ is a Hilbert space with the norm
\[
\|f\|_{\cD((-\tri)^\ga)}:=\left(\sum_{n=1}^\infty\left|\la_n^\ga(f,\vp_n)\right|^2\right)^{1/2},\quad f\in\cD((-\tri)^\ga).
\]

For $1\le p\le\infty$ and a Banach space $X$, we say that $f\in L^p(0,T;X)$ if
\[
\|f\|_{L^p(0,T;X)}:=\left\{\!\begin{alignedat}{2}
& \left(\int_0^T\|f(\,\cdot\,,t)\|_X^p\,\rd t\right)^{1/p} & \quad & \mbox{if }1\le p<\infty\\
& \mathop{\mathrm{ess}\sup}_{0<t<T}\|f(\,\cdot\,,t)\|_X & \quad & \mbox{if }p=\infty
\end{alignedat}\right\}<\infty.
\]

We define the forward Riemann-Liouville integral operator of order $\al\in(0,1)$:
\[
J_{0+}^\al f(t):=\f1{\Ga(\al)}\int_0^t\f{f(s)}{(t-s)^{1-\al}}\,\rd s.
\]
Then it is easily seen that the Caputo derivative $\pa_t^\al f=J_{0+}^{1-\al}(\f{\rd f}{\rd t})$. For the solution representation, we define the Mittag-Leffler function (e.g., Podlubny \cite{P99})
\[
E_{\al,\be}(z):=\sum_{k=0}^\infty\f{z^k}{\Ga(\al k+\be)},\quad z\in\BC,\ \al>0,\ \be>0.
\]
The following estimate is later useful:
\begin{equation}\label{eq-est-ML}
|E_{\al,\be}(-\eta)|\le\f C{1+\eta},\quad\eta\ge0,\ 0<\al<2,\ \be>0,
\end{equation}
where $C>0$ is a constant depending only on $\al$ and $\be$.

For inverse problems, we have to study the unique existence of solution $u(x,t)$ to the initial-boundary value problem \eqref{eq-ibvp-u} and its properties, which is called the forward problem, contrasted with the inverse problems.

Now we collect the basic results on the forward problem
\begin{equation}\label{eq-ibvp-w}
\begin{cases}
(\pa_t^\al-\tri)w(x,t)=F(x,t), & x\in\Om,\ 0<t<T,\\
w(x,0)=a(x), & x\in\Om,\\
w(x,t)=0, & x\in\pa\Om,\ 0<t<T.
\end{cases}
\end{equation}

\begin{lemma}\label{lem-FP-u}
{\rm(a)} Let $F=0$ and $a\in\cD((-\tri)^\ga)$ with $\ga\ge0$. Then there exists a unique solution $w\in\bigcap_{1\le p\le\infty} L^p(0,T;\cD((-\tri)^{\ga+\f1p}))$ to \eqref{eq-ibvp-w}, where we interpret
$\f1p=0$ for $p=\infty$. Moreover, the solution $w$ allows the representation
\[
w(\,\cdot\,,t)=\sum_{n=1}^\infty E_{\al,1}(-\la_nt^\al)(a,\vp_n)\vp_n
\]
in $L^p(0,T;\cD((-\tri)^{\ga+\f1p}))$ for any $p\in[1,\infty]$, and $w:(0,T]\longrightarrow L^2(\Om)$ can be analytically extended to a sector $\{z\in\BC;\ z\ne0,\ |\mathrm{arg}\,z|<\f\pi2\}$. Furthermore, there exists a constant $C>0$ such that for any $t>0$,
\begin{align*}
& \|w(\,\cdot\,,t)\|_{\cD((-\tri)^{\ga+1/p})}\le C\,t^{-\f\al p}\|a\|_{\cD((-\tri)^\ga)},\\
& \|w(\,\cdot\,,t)\|_{\cD((-\tri)^{\ga+1})}+\|\pa_t^\al w(\,\cdot\,,t)\|_{\cD((-\tri)^\ga)}\le C\,t^{-\al}\|a\|_{\cD((-\tri)^\ga)}.
\end{align*}

{\rm(b)} Let $a=0$ and $F\in L^p(0,T;\cD((-\tri)^\ga))$ with $p\in[1,\infty]$ and $\ga\ge0$. Then there exists a unique solution $w\in\bigcap_{0<\ve\le1}L^p(0,T;\cD((-\tri)^{\ga+1-\ve}))$ to \eqref{eq-ibvp-w}. Moreover, the solution $w$ allows the representation
\begin{equation}\label{eq-rep-u}
w(\,\cdot\,,t)=\sum_{n=1}^\infty\left(\int_0^t s^{\al-1}E_{\al,\al}(-\la_n s^\al)\,(F(\,\cdot\,,t-s),\vp_n)\,\rd s\right)\vp_n
\end{equation}
in $L^p(0,T;\cD((-\tri)^{\ga+1-\ve}))$ for any $\ve\in(0,1]$. Furthermore, there exists a constant $C>0$ such that for $\forall\,\ve\in(0,1]$,
\[
\|w\|_{L^p(0,T;\cD((-\tri)^{\ga+1-\ve}))}+\|\pa_t^\al w\|_{L^p(0,T;\cD((-\tri)^{\ga-\ve}))}\le\f C\ve\|F\|_{L^p(0,T;\cD((-\tri)^\ga))}.
\]
\end{lemma}

The above well-posedness results are refinements of that stated in \cite[Theorems 2.1--2.2]{SY11a}, which can be easily verified by the arguments in Li, Liu and Yamamoto \cite{LLY15}.

Especially, Lemma \ref{lem-FP-u} reflects the limited smoothing property of time-fractional diffusion equations, that is, in the case of $F\equiv0$, the improvement of the spatial regularity of solution is at most $2$, compared with the initial value. More precisely, the regularity improvement of the homogeneous problem can exactly reach $2$ at the cost of a weakened norm in time. On the other hand, with a source term of $L^p$ regularity in time, the regularity improvement of the inhomogeneous one can never reach $2$ except for the special case of $p=2$, where the complete monotonicity property of Mittag-Leffler functions can be utilized.

Henceforth, we understand the class of solutions as described in Lemma \ref{lem-FP-u}.

We can refer to Eidelman and Kochubei \cite{EK04} on fundamental solutions, Gorenflo, Luchko and Yamamoto \cite{GLY15} on the suitable function spaces for solutions, Kubica and Yamamoto \cite{KY17} on the weak solution and improved regularity, Luchko \cite{L09} on a solution formula and the maximum principle, Sakamoto and Yamamoto \cite{SY11a} on the well-posedness including some inverse problems by the representation of solutions, Zacher \cite{Z09} on generalized treatments on weak solutions. Here we represent very limited references and we can consult also other related chapters of this handbook.

We can represent the solution to \eqref{eq-ibvp-u} by a fractional Duhamel's principle (see Liu, Rundell and Yamamoto \cite{LRY16}).

\begin{lemma}\label{lem-Duhamel}
Let $u$ be the solution to \eqref{eq-ibvp-u}, where $\rho\in L^1(0,T)$ and $g\in\cD((-\tri)^\ga)$ with $\ga\ge0$. Then $u$ allows the representation
\begin{equation}\label{eq-Duhamel}
J_{0+}^{1-\al}u(\,\cdot\,,t)=\int_0^t\rho(t-s)v(\,\cdot\,,s)\,\rd s,
\end{equation}
where $J_{0+}^{1-\al}$ is the $(1-\al)$-th order Riemann-Liouville integral, and $v$ solves the homogeneous problem
\begin{equation}\label{eq-ibvp-v}
\begin{cases}
(\pa_t^\al-\tri)v(x,t)=0m & x\in\Om,\ t>0,\\
v(x,0)=g(x), & x\in\Om,\\
v(x,t)=0, & x\in\pa\Om,\ t>0.
\end{cases}
\end{equation}
\end{lemma}

Lemma \ref{lem-Duhamel} relates the inhomogeneous problem \eqref{eq-ibvp-u} with the homogeneous one \eqref{eq-ibvp-v}.  Therefore, it suffices to study \eqref{eq-ibvp-v} for the inverse source problems for
\eqref{eq-ibvp-u}. Later we will see that Lemma \ref{lem-Duhamel} acts as the starting point for discussing the inverse problems. For the fractional Duhamel's principle, see also Umarov and Saidamatov \cite{US07}, Zhang and Xu \cite{ZX11}.

\section{Inverse $t$-source problems}\label{sec-ISP-t}

\subsection{Two-sided Lipschitz stability in the case of $x_0\in\supp\,g$}

\begin{theorem}\label{thm-SY11a}
Let $g\in\cD((-\tri)^\ga)$ with $\ga>\f d4$ and $u$ satisfy \eqref{eq-ibvp-u} for $\rho\in C[0,T]$. If $g(x_0)\ne0$, then there exists a constant $C>0$ such that
\begin{equation}\label{eq-SY11a-Lip}
C^{-1}\|\pa_t^\al u(x_0,\,\cdot\,)\|_{L^\infty(0,T)}\le\|\rho\|_{C[0,T]}\le C\|\pa_t^\al u(x_0,\,\cdot\,)\|_{L^\infty(0,T)}.
\end{equation}
\end{theorem}

By the Sobolev embedding $\cD((-\tri)^\ga)\subset C(\ov\Om)$ with $\ga>\f d4$,the regularity of the known spatial component $g$ in the above theorem turns out to be sufficient for a pointwise definition of $g(x)$. Theorem \ref{thm-SY11a} slightly improves the regularity $g\in\cD((-\tri)^{\ga_0})$ with $\ga_0>\f{3d}4+1$ in \cite[Theorem 4.4]{SY11a}. Here we briefly explain how to realize such a reduction in regularity.

Actually, the first inequality in \eqref{eq-SY11a-Lip} is a direct corollary of Lemma \ref{lem-FP-u}(b) and the Sobolev embedding theorem. In order to show the second inequality in \eqref{eq-SY11a-Lip}, we substitute the representation
\eqref{eq-rep-u} into the governing equation in \eqref{eq-ibvp-u} and substitute $x=x_0$ to write formally
\begin{equation}\label{eq-ISP-x0}
\rho(t)=\f1{g(x_0)}\left\{\pa_t^\al u(x_0,t)+\int_0^t Q(x_0,s)\rho(t-s)\,\rd s\right\},
\end{equation}
where
\[
Q(\,\cdot\,,t):=t^{\al-1}\sum_{n=1}^\infty\la_n E_{\al,\al}(-\la_nt^\al)(g,\vp_n)\vp_n.
\]
Introducing $\ve:=\f12(\ga-\f d4)>0$, we employ \eqref{eq-est-ML} to estimate
\begin{align*}
\|Q(\,\cdot\,,t)\|_{\cD((-\tri)^{d/4+\ve})}^2 & =t^{2(\al-1)}\sum_{n=1}^\infty\left|\la_n^{1-\ve}E_{\al,\al}(-\la_n t^\al)\right|^2\left|\la_n^{2\ve+d/4}(g,\vp_n)\right|^2\\
& \le C^2\,t^{2(\al-1)}\sum_{n=1}^\infty\left(\f{(\la_n t^\al)^{1-\ve}}{1+\la_n t^\al}t^{\al(\ve-1)}\right)^2|\la_n^\ga(g,\vp_n)|^2\\
& \le\left(C\|g\|_{\cD((-\tri)^\ga)}t^{\al\ve-1}\right)^2.
\end{align*}
Using the Sobolev embedding $\cD((-\tri)^{\f d4+\ve})\subset H^{d/2+2\ve}(\Om)\subset C(\ov\Om)$, we obtain
\[
|Q(x_0,t)|\le\|Q(\,\cdot\,,t)\|_{C(\ov\Om)}\le C\|Q(\,\cdot\,,t)\|_{\cD((-\tri)^{d/4+\ve})}\le C\|g\|_{\cD((-\tri)^\ga)}t^{\al\ve-1}.
\]
Therefore, the above estimate implies
\[
|\rho(t)|\le C\|\pa_t^\al u(x_0,\,\cdot\,)\|_{L^\infty(0,T)}+C\|g\|_{\cD((-\tri)^\ga)}\int_0^t s^{\al\ve-1}|\rho(t-s)|\,\rd s.
\]
Applying a Gronwall-type inequality in \cite[Lemma 7.1.1]{H81}, we complete the proof of the second inequality in \eqref{eq-SY11a-Lip}.\medskip

In the same direction, also several other papers obtained Lipschitz stability in slightly different formulations. By assuming the homogeneous Neumann condition instead of that in \eqref{eq-ibvp-u}, the observation point $x_0$
can be placed on the boundary, and we have the following result.

\begin{theorem}[Wei, Li and Li \cite{WLL16}]
Let $u$ be the solution to \eqref{eq-ibvp-u} with the homogeneous Neumann boundary condition, where $\rho$ is absolutely continuous on $[0,T]$ and $g\in\cD((-\tri+1)^\ga)$ with $\ga>\f d2+1$. If $x_0\in\pa\Om$ and $g(x_0)\ne0$, then \eqref{eq-SY11a-Lip} still holds.
\end{theorem}

By a similar argument as that of Theorem \ref{thm-SY11a}, we can also reduce the regularity assumption in the above theorem, and we skip the details here.

On the other hand, the following result reveals that Theorem \ref{thm-SY11a} holds true with a more general source term.

\begin{theorem}[Fujishiro and Kian \cite{FK16}]\label{thm-FK16}
Let $u$ be the solution to \eqref{eq-ibvp-u} with $F(x,t)=\rho(t)g(x,t)$, where $g\in L^p(0,T;\cD(-\tri))$ with $\f8\al<p\le\infty$ and there exists a constant $\de>0$ such that $|g(x_0,\,\cdot\,)|\ge\de$ a.e.\! in $(0,T)$. Then there exists a constant $C>0$ such that
\[
\|\rho\|_{L^p(0,T)}\le C\|\pa_t^\al u(x_0,\,\cdot\,)\|_{L^p(0,T)}.
\]
\end{theorem}

The key to the proof is an $L^p$ estimate for \eqref{eq-ISP-x0} (see \cite[Lemma 4]{FK16}), and here we omit the details. In \cite{FK16}, Theorem \ref{thm-FK16} is applied for establishing the conditional stability for a corresponding inverse coefficient problem by the same observation data. For $1<\al<2$, Wu and Wu \cite{WW14b} obtained the uniqueness in a more general formulation under the same non-vanishing condition.\medskip

Ruan and Wang \cite{RW17} adopts distributed observations: given $\si\in C_0^\infty(\Om)$ satisfying $\si\ge0$ and
$\si\not\equiv0$, one measures
\[
\int_\Om u(x,t)\si(x)\,\rd x,\quad0<t<T.
\]
Under the condition
\[
(\si,g)=\int_\Om\si(x)g(x)\,\rd x\ne0,
\]
\cite[Theorem 1]{RW17} established a similar estimate to \eqref{eq-SY11a-Lip} in fractional Sobolev norms. In this direction, see also Aleroev, Kirane and Malik \cite{AKM13}, which restricts $\si\equiv1$ but generalizes $F(x,t)=\rho(t)g(x,t)$ in \eqref{eq-ibvp-u}.

\subsection{Uniqueness and stability with general observation point $x_0$}

First we state a uniqueness result for any $x_0\in\Om$ not necessarily satisfying $g(x_0)\ne0$, which removes the restriction on the space dimensions and reduces the required regularity of $\rho$ in the result in Liu, Rundell and Yamamoto \cite{LRY16}.

\begin{theorem}\label{thm-ISPt-unique}
We assume that $\rho\in L^1(0,T)$, $g\in\cD((-\tri)^\ga)$ with $\ga>\f d4-1$, $g\ge0$ and $g\not\equiv0$. Then
\[
u(x_0,t)=0,\quad0<t<T\quad\mbox{implies}\quad\rho(t)=0,\quad0<t<T.
\]
\end{theorem}

The key to such an improvement in Theorem \ref{thm-ISPt-unique} is Lemma \ref{lem-Duhamel}, that is, a weak form of the fractional Duhamel's principle. Taking $x=x_0$ in \eqref{eq-Duhamel} of Lemma \ref{lem-Duhamel} and using $u(x_0,\,\cdot\,)=0$ in $(0,T)$, we have
\[
0=J_{0+}^{1-\al}u(x_0,t)=\int_0^t\rho(t-s)v(x_0,s)\,\rd s,\quad0<t<T.
\]
After this step, one can follow the arguments in \cite{LRY16} to employ the Titchmarsh convolution theorem (see Titchmarsh \cite{T26}) and some strict positivity property of the solution to \eqref{eq-ibvp-v} (see \cite{LRY16}) to conclude the result.

Moreover Liu \cite{L17} generalized Theorem \ref{thm-ISPt-unique} for multi-term time-fractional diffusion equations.\medskip

Next, we continue to investigate the stability of the inverse $t$-source problem especially in the case of $x_0\not\in\supp\,g$. Only in this part, instead of the initial-boundary value problem \eqref{eq-ibvp-u}, we consider the Cauchy problem in the whole space:
\begin{equation}\label{eq-Cauchy-u}
\begin{cases}
(\pa_t^\al-\tri)u(x,t)=g(x)\rho(t), & x\in\BR^d,\ 0<t<T,\\
u(x,0)=0, & x\in\BR^d.
\end{cases}
\end{equation}

In order to state the stability, for given constants $M>0$ and $N\in\{0,1,2,\ldots\}$, we define the admissible set of the unknown temporal components $\rho$ by
\begin{align}
\cU_N:=\{f\in C^1[0,T];\; & \|f\|_{C^1[0,T]}\le M,\nonumber\\
& f\mbox{ changes signs at most }N\mbox{ times on }(0,T)\}.\label{eq-def-sign}
\end{align}
This admissible set was introduced in Saitoh, Tuan and Yamamoto \cite{STY02}, where the same inverse $t$-source problem is discussed for $\al=1$.

With the above preparations, now we can state the main stability result.

\begin{theorem}[Liu and Zhang \cite{LZ17}]\label{thm-ISPt-stability}
Let $u$ satisfy \eqref{eq-Cauchy-u}, where we assume that
\begin{equation}\label{eq-asp-g}
g\ge0,\quad g\not\equiv0,\quad\begin{cases}
g:\mbox{bounded continuous} & \mbox{if }d=1,\\
g:\mbox{locally H\"older continuous} &  \mbox{if }d=2,3.
\end{cases}
\end{equation}
Let $\cU_N$ be defined by \eqref{eq-def-sign} with given constants $M>0$ and $N\in\{0,1,2,\ldots\}$.

{\rm(a)} Let $N=0$. Then for any $\de_1\in(0,T)$, there exists a constant $B_{\de_1}>0$ depending only on $x_0,g$ such that $B_{\de_1}\to\infty$ as $\de_1\to0$ and
\begin{equation}\label{eq-stab-a}
\|\rho\|_{L^1(0,T-\de_1)}\le\f{T^{1-\al}B_{\de_1}}{\Ga(2-\al)}\|u_\rho(x_0,\,\cdot\,)\|_{L^1(0,T)}
\end{equation}
for all $\rho\in\cU_0$.

{\rm(b)} Let $N\in\{1,2,3,\ldots\}$. We further assume that $x_0\not\in\supp\,g$ and $\|u(x_0,\,\cdot\,)\|_{L^1(0,T)}$ is sufficiently small. Then there exist a sufficiently small $\de_2>0$ and a constant $C>0$ such that
\begin{equation}\label{eq-stab-d}
\|\rho\|_{L^1(0,T-\de_2)}\le C\,(\underbrace{\log(\cdots(\log|\log}_{N+1}\|u_\rho(x_0,\,\cdot\,)\|_{L^1(0,T)}|)\cdots))^{-2(2/\al-1)}
\end{equation}
for all $\rho\in\cU_N$.
\end{theorem}

In the conditional stability \eqref{eq-stab-a} and \eqref{eq-stab-d}, we need to take data over a longer observation time interval $(0,T)$ for estimating $\|\rho\|_{L^1(0,T-\delta_k)}$ with $k=1,2$. Our stability estimate is weaker when the time $N$ of changing signs is increasing.

The parameter $B_{\de_1}>0$ in Theorem \ref{thm-ISPt-stability}(a) turns out to be
\[
B_{\de_1}:=\f1{\|v(x_0,\,\cdot\,)\|_{L^1(0,\de_1)}}=\left(\int_0^{\de_1} v(x_0,t)\,\rd t\right)^{-1},
\]
where $v$ solves the homogeneous Cauchy problem
\begin{equation}\label{15}
\begin{cases}
(\pa_t^\al-\tri)v(x,t)=0, & x\in\BR^d,\ t>0,\\
v(x,0)=g(x)\ge0,\not\equiv0, & x\in\BR^d.
\end{cases}
\end{equation}
For the details, we refer to \cite[Section 2]{LZ17}.

According to Eidelman and Kochubei \cite{EK04}, $v$ is strictly positive at least for sufficiently small $t>0$. Furthermore \cite{EK04} gives a lower bound for $v$, which validates the quantitative analysis when $\rho$ does change signs.\medskip

More precisely, the proof of Theorem \ref{thm-ISPt-stability} is based on the following lemma which one can prove also by using \cite{EK04}.

\begin{lemma}\label{lem-Cauchy}
{\rm(a)} Let $u$ be the solution to \eqref{eq-Cauchy-u}, where we assume $\rho\in C[0,\infty)$ and \eqref{eq-asp-g}. Then Lemma $\ref{lem-Duhamel}$ still holds, that is, $u$ allows the same representation \eqref{eq-Duhamel}, where $v$ solves
the initial value problem \eqref{15} for the homogeneous equation.

{\rm(b)} Let $g$ satisfy \eqref{eq-asp-g}. Then there exists a classical solution to \eqref{15}, which takes the form
\[
v(x,t)=\int_{\BR^d}K_\al(x-\xi,t)g(\xi)\,\rd\xi,
\]
where the fundamental solution $K_\al(x,t)$ satisfies the following asymptotic behavior as $t\downarrow0$: If $|x|>r$ for some fixed $r>0$, then there exist a constant $C>0$ depending on $\al,d,r$ such that
\[
K_\al(x,t)\sim t^{-\f{\al d}{2(2-\al)}}|x|^{-\f{d(1-\al)}{2-\al}}\exp\left(-C\,t^{-\f\al{2-\al}}|x|^{\f2{2-\al}}\right).
\]
\end{lemma}

For other related works on the inverse $t$-source problem, see also Jin and Rundell \cite{JR15}, Wang and Wu \cite{Wang-Wu14}. At the end of this section, we briefly mention the numerical reconstruction method for the inverse $t$-source problem developed in \cite{LZ17}. In order to specify the dependency on $\rho$, by $u(\rho)$ we denote the solution of \eqref{eq-ibvp-u}, and let $\rho_*$ be the true solution. Under the same non-negativity assumption of $g$ as before, we propose the fixed-point iteration
\[
\rho_m(t)=\left\{\!\begin{alignedat}{2}
& 0, & \quad & m=0,\\
& \rho_{m-1}(t)+\f{\pa_t^\al(u(\rho_*)-u(\rho_{m-1}))(x_0,t)}K, & \quad &
m=1,2,\ldots,
\end{alignedat}\right.
\]
where $K>0$ is a constant such that $\|v(x_0,\,\cdot\,)\|_{C[0,T]}\le K$. Since $v$ can be computed in advance, we can easily evaluate $K$, and the proposed iteration only involves one-dimensional computation in time by taking
time derivatives in \eqref{eq-Duhamel} of Lemma \ref{lem-Duhamel}. See \cite{LZ17} for further details on the convergence analysis and numerical examples.

\section{Inverse $x$-source problems}\label{sec-ISP-x}

In this section, we investigate inverse $x$-source problems for \eqref{eq-ibvp-u}.

\subsection{Final observation data}

For the inverse $x$-source problem with final time observation data, we review three relevant results. First, we introduce the conditional H\"older stability obtained in Wang, Zhou and Wei \cite{WZW13}:

\begin{theorem}\label{thm-Wei}
Let $u$ be the solution to \eqref{eq-ibvp-u}, where we assume $g\in\cD((-\tri)^\ga)$ with some fixed $\ga>0$ and $\rho(t)\equiv1$. If $g$ satisfies an a priori estimate $\|g\|_{\cD((-\tri)^\ga)}\le E$ with a constant $E>0$, then there exists a constant $C>0$ such that
\[
\|g\|_{L^2(\Om)}\le C\,E^{\f1{\ga+1}}\|u(\,\cdot\,,T)\|_{L^2(\Om)}^{\f\ga{\ga+1}}.
\]
\end{theorem}

Similarly to several previous theorems, the above theorem follows immediately from Lemmata \ref{lem-FP-u} and \ref{lem-Duhamel} and the key estimate \eqref{eq-est-ML} in Section \ref{sec-pre}. Such results as Theorem \ref{thm-Wei} are called to be conditional stability, because in order to estimate the $L^2$ norm of $g$, one should assume its a priori bound with some norm.\medskip

In general, it is technically difficult to establish unconditional stability for the inverse $x$-source problem for general $\al\in(0,1)$. Only in the case of $\al=\f12$, Yamamoto and Zhang \cite{YZ12} proved such stability as a by-product for treating a corresponding inverse coefficient problem. For conciseness, here we only discuss a slightly simpler formulation
than that in \cite{YZ12}:
\begin{equation}\label{eq-gov-CE}
\begin{cases}
(\pa_t^{\f12}-\pa_x^2)u(x,t)=g(x)\rho(x,t), & 0<x<\ell,\ 0<t<T,\\
u(x,0)=0, & 0<x<\ell,\\
u(0,t)=\pa_xu(0,t)=0, & 0<t<T.
\end{cases}
\end{equation}

We state the local H\"older stability for the inverse $x$-source problem in this case.

\begin{theorem}\label{thm-half}
Let $u$ satisfy \eqref{eq-gov-CE} and suppose that $u$, $g$ and $\rho$ are sufficiently smooth. Fix $t_0\in(0,T)$ arbitrarily and choose $x_0>\ell$ such that $x_0-\ell>0$ is sufficiently small. Define the level set $\Om_\ve:=\{x\in(0,\ell);\,|x-x_0|^2>\ve\}$ for $\ve>0$. If $\rho(\,\cdot\,,t_0)\ne0$ in $\ov\Om$ and $g(0)=g'(0)=0$, then there exist constants $C>0$ and $\ka\in(0,1)$ such that
\[
\|g\|_{H^2(\Om_{3\ve})}\le C\|u(\,\cdot\,,t_0)\|_{H^4(\Om_\ve)}^\ka.
\]
\end{theorem}

The key idea to prove the above theorem is a transform of \eqref{eq-gov-CE} to an equation governed by the fourth-order differential operator $\pa_t-\pa_x^4$. For such an equation, one can apply a class of weighted $L^2$ estimates called Carleman estimates to prove the stability of inverse problems; see also Xu, Cheng and Yamamoto \cite{XCY11} for the derivation of the Carleman estimate. Unfortunately, the idea of transforming to an equation of integer order works at most for $\al\in\mathbb Q$. Actually, even the case of $\al=\f13$ involves huge amounts of calculations in applying Carleman estimates. Also in Chapter ``Inverse problems of determining coefficients of the fractional partial differential equations'' of this handbook, we describe about Carleman estimates for fractional diffusion equations in general spatial dimensions. Here we refer to Cheng, Lin and Nakamura \cite{CLN13} for $\al=\f12$, Lin and Nakamura \cite{LN16} for $\al\in(0,1)$, and Lin and Nakamura \cite{LN2018} for equations with multi-term time fractional detivatives of orders $\al\in(0,2)\setminus\{1\}$. Their Carleman estimates yield the unique continuation for Caputo time-fractional diffusion equations, but are not applicable to inverse problems.\medskip

Next, we discuss the structure of the inverse problem with final observation data by the Fredholm alternative. We consider a more general inverse source problem when $\rho$ depends also on $x$:
\[
\begin{cases}
(\pa_t^\al-\tri-p(x))u(x,t)=g(x)\rho(x,t), & x\in\Om,\ 0<t<T,\\
u(x,0)=0, & x\in\Om,\\
u(x,t)=0, & x\in\pa\Om,\ 0<t<T.
\end{cases}
\]

With suitable regularity assumption on $\rho$, we reduce this inverse source problem with final data to a Fredholm equation of the second kind:
\[
g+Kg=-\f{(\tri+p(x))u(x,T)}{\rho(x,T)},\quad x\in\Om,
\]
where $K:L^2(\Om)\longrightarrow L^2(\Om)$ is a compact operator (e.g., Sakamoto and Yamamoto \cite{SY11b}). Therefore, if $-1$ is not an eigenvalue of $K$, then the inverse source problem is well-posed in the sense of Hadamard. In the case where $-1$ is an eigenvalue of $K$, the non-uniqueness for the inverse source problem is restricted only in a finite dimensional subspace of $L^2(\Om)$. This property by the Fredholm alternative was originally proved for the inverse problem with the final observation for parabolic equations (e.g., \cite{POV00}) and the same property holds also for the fractional diffusion equation (e.g., \cite{SY11b}).\medskip

Finally, we treat the problem from a different point of view. Let us consider the perturbation of the governing equation in
\eqref{eq-ibvp-u} with a parameter $r>0$:
\begin{equation}\label{eq-perturb}
\begin{cases}
\pa_t^\al u(x,t)=r^\al(\tri+p(x))u(x,t)+g(x)\rho(x,t), & x\in\Om,\ 0<t<T,\\
u(x,0)=0, & x\in\Om,\\
u(x,t)=0, & x\in\pa\Om,\ 0<t<T.
\end{cases}
\end{equation}
Specifying the dependency on $r$ and given smooth $g$, by $u(r,g)$ we denote the unique solution to \eqref{eq-perturb}. We state the generic well-posedness result in the following theorem, which is a slight simplification of the main theorem
in Sakamoto and Yamamoto \cite{SY11b}.

\begin{theorem}\label{thm-generic}
Let $u(r,g)$ be the solution to \eqref{eq-perturb}. We assume that $p\in C(\ov\Om)$, $\rho\in C^1([0,T];L^\infty(\Om))$ and $\rho(\,\cdot\,,T)\ne0$. For any open interval $I\in(0,\infty)$, there exists a finite set $E=E(\al,\rho,I)\subset I$ such that for any $r\in I\setminus E$ and $\psi\in\cD(-\tri)$, there exists a unique solution $\{u(r,g),g\}\in C([0,T];\cD(-\tri))\times L^2(\Om)$ to \eqref{eq-perturb} satisfying $u(r,g)(\,\cdot\,,T)=\psi$. Moreover, there exists a constant $C>0$ such that
\[
\|g\|_{L^2(\Om)}+\|u(r,g)\|_{C([0,T];H^2(\Om))}+\|\pa_t^\al u(r,g)\|_{C([0,T];L^2(\Om))}\le C\|\psi\|_{H^2(\Om)}.
\]
\end{theorem}

We do not know whether $E=\emptyset$. We can understand Theorem \ref{thm-generic} as follows. For an arbitrarily given target function $\psi\in\cD(-\tri)$, we attempt to find a pair $(r,g)\in(0,\infty)\times L^2(\Om)$ such that $u(r,g)(\,\cdot\,,T)=\psi$. Unfortunately, with an arbitrarily fixed $r>0$, for example, $r=1$, we do not know whether there exists a unique $g\in L^2(\Om)$ satisfying the above condition. However, by taking $I$ as any open neighborhood of $r=1$, Theorem \ref{thm-generic} asserts that the problem may be ill-posed only for a finite set of $I$. This inverse problem is generically well-posed in the sense of Hadamard, although we do not know whether the original problem with $r=1$ is well-posed. We refer to Choulli and Yamamoto \cite{CY96,CY99} for the generic well-posedness for inverse parabolic problems.

With the aid of the Fredholm alternative, the key to proving Theorem \ref{thm-generic} is the analytic perturbation
theory (see Kato \cite{K76}).

In this direction, Tatar and Ulusoy \cite{TU14} discusses the same type of inverse $x$-source problem with final observation data for
\[
(\pa_t^\al+r^\al(-\tri)^{\be/2})u(x,t)=g(x)\rho(t),\quad x\in\Om,\ 0<t<T
\]
with given $\be\in(0,2)$, and see Tatar, Tinaztepe and Ulusoy \cite{TTU2016} as for a numerical method.

For other related works on the inverse $x$-source problem with final observation, see also Kawamoto \cite{K18}, Kirane and Malik \cite{KM11}, Kirane, Malik and Al-Gwaiz \cite{KMA13}.

\subsection{Uniqueness by partial interior observation}

Now we continue to study the inverse $x$-source problem with interior observation data. Regarding the same type of problems for classical partial differential equations such as wave and heat equations, we know that there are quite a lot of stability results based on Carleman estimates (e.g., Bellassoued and Yamamoto \cite{BY17}, Klibanov and Timonov \cite{KT04}, Yamamoto \cite{Y09}). However, except for rather special cases like $\al=\f12$, such methodology does not work for fractional equations due to the absence of convenient formulae of integration by parts for fractional derivatives. As a result, to the best of our knowledge, the stability of the inverse $x$-source problem with general $\al\in(0,1)$ mostly keeps open, and here we can only review the uniqueness result in Jiang, Li, Liu and Yamamoto \cite{JLLY17}.

\begin{theorem}[Jiang, Li, Liu and Yamamoto \cite{JLLY17}]\label{thm-ISPx-unique}
Let $g\in L^2(\Om)$ and assume that $\rho\in C^1[0,T]$ with $\rho(0)\ne0$. Let $u$ be the solution to \eqref{eq-ibvp-u} and $\om\subset\Om$ be an arbitrary nonempty subdomain. Then
\[
u=0\mbox{ in }\om\times(0,T)\quad\mbox{implies}\quad g=0\mbox{ in }\Om.
\]
\end{theorem}

The keys to proving the above theorem are the Duhamel's principle \eqref{eq-Duhamel} and the following uniqueness for \eqref{eq-ibvp-v} in \cite{JLLY17}:

\begin{lemma}\label{lem-UCP}
Let $v$ be the solution to \eqref{eq-ibvp-v} with $g\in L^2(\Om)$, and $\om\subset\Om$ be an arbitrary nonempty subdomain. Then
\[
v=0\mbox{ in }\om\times(0,T)\quad\mbox{implies}\quad v=0\mbox{ in }\Om\times(0,\infty).
\]
\end{lemma}

We briefly introduce a numerical method for the inverse $x$-source problem developed in \cite{JLLY17}. In order to specify the dependency on $g$, let us denote the solution of \eqref{eq-ibvp-u} by $u(g)$, and let $g_*$ be the true solution. Then we propose the iterative thresholding algorithm
\[
g_{m+1}=\f K{K+\be}g_m-\f1{K+\be}\int_0^T\rho\,z(g_m)\,\rd t,\quad m=0,1,2,\ldots,
\]
where $K>0$ and $\be>0$ are suitably chosen parameters. Here $z(g_m)$ solves the backward problem
\[
\begin{cases}
-J_{T-}^{1-\al}(\pa_t z)-\tri z=\chi_\om\left(u(g_m)-u(g_*)\right) & \mbox{in }\Om\times(0,T),\\
z=0 & \mbox{in }\Om\times\{T\},\\
z=0 & \mbox{on }\pa\Om\times(0,T),
\end{cases}
\]
where $\chi_\om$ is the characteristic function of $\om$ and $J_{T-}^{1-\al}$ denotes the backward Riemann-Liouville integral operator defined by
\[
J_{T-}^\al f(t):=\f1{\Ga(\al)}\int_t^T\f{f(s)}{(s-t)^{1-\al}}\,\rd s.
\]
See \cite{JLLY17} for further details.\medskip

We mention that there are other kinds of observation data. For example, Zhang and Xu \cite{ZX11} discussed the determination of $g$ in
\[
\begin{cases}
(\pa_t^\al-\pa_x^2)u(x,t)=g(x), & 0<x<1,\ 0<t<T,\\
u(x,0)=0, & 0<x<1,\\
\pa_x u(0,t)=\pa_x u(1,t)=0, & 0<t<T
\end{cases}
\]
by the boundary observation $u(0,t)$. The uniqueness can be easily shown by using Lemmata \ref{lem-FP-u} and \ref{lem-Duhamel} and Laplace transform. The key is the $t$-analyticity of the solution $u(x,t)$, which is guaranteed by $\rho\equiv1$, and their argument does not work for general $\rho\in C[0,T]$.

In higher spatial dimensions, Wei, Sun and Li \cite{WSL16} studied an inverse $x$-source problem by extra boundary data for $0<t<\infty$ and establishes the uniqueness in the inverse problem. In the case where we consider the initial-boundary value problem and data over the infinite time interval $0<t<\infty$, similarly to \cite{ZX11} we can take the Laplace transforms, so that the uniqueness follows.

\section{Related inverse source problems}\label{sec-re}

In inverse problems, we are required to determine various quantities such as source terms, coefficients and parameters describing the fractional derivatives and so there are naturally various types of inverse problems. Thus it is not reasonable to rigorously classify the inverse problems for fractional partial differential equations and in this section, we supplementarily survey inverse problems of determining boundary quantities, some of which is the determination of source term located on the boundary $\pa\Om$.

\subsection{Inverse problem of determining $\rho(t)$ in the boundary source}

Let $\ga\subset\pa\Om$ be a relatively open subboundary. We consider
\[
\begin{cases}
\pa_t^\al u(x,t)=\tri u(x,t), & x\in\Om,\ 0<t<T,\\
u(x,t)=\begin{cases}
\rho(t)g(x),\\
0,
\end{cases} & \!\begin{aligned}
& x\in\ga,\ 0<t<T,\\
& x\in\pa\Om\setminus\ov\ga,\ 0<t<T,
\end{aligned}\\
u(x,0)=0, & x\in\Om.
\end{cases}
\]
Then, given $g$, determine $\rho(t)$, $0<t<T$ by distributed measurement data
\[
\int_\Om u(x,t)\si(x)\,\rd x,\quad0<t<T.
\]
Here $\si\in C^\infty_0(\Om)$, $\ge0$, $\not\equiv0$ in $\Om$, is an arbitrarily chosen function, which describes a weight. As an extremal case, setting $\si(x)=\de(x-x_0)$ with fixed $x_0\in\Om$, we can reduce our data to the pointwise $u(x_0,t)$, $0<t<T$. However, such pointwise data are not studied.

Liu, Yamamoto and Yan \cite{LYY16} proves

\begin{theorem}\label{thm9}
Let $\f12<\al<1$ and $g\in C^\infty_0(\ga)$ satisfy $g\ge0$, $\not\equiv0$ on $\ga$. Then there exists a constant $C>0$ such that
\[
C^{-1}\left\|\int_\Om u(x,\,\cdot\,)\si(x)\,\rd x\right\|_{H^\al(0,T)}\le\|\rho\|_{H^\al(0,T)}\le C\left\|\int_\Om u(x,\,\cdot\,)\si(x)\,\rd x\right\|_{H^\al(0,T)}
\]
for all $\rho\in H^\al(0,T)$ satisfying $\rho(0)=0$.
\end{theorem}

Here $H^\al(0,T)$, $0<\al<1$, is the fractional Sobolev space defined by Sobolev-Slobodecki norm (e.g., Adams \cite{A75}). By the Sobolev embedding, we know that $H^\al(0,T)\subset C[0,T]$ for $\al>\f12$, and so $\rho(0)=0$ makes sense for $\rho\in H^\al(0,T)$ with $\al>\f12$. The paper \cite{LYY16} considers a more general elliptic operator instead of $-\tri$, and gives numerical methods for reconstructing $\rho(t)$ by noisy data for $\int_\Om u(x,t)\si(x)\,\rd x$ for $0<t<T$.

\subsection{Determination of $t$-coefficient in the Robin boundary condition}

For
\[
\begin{cases}
\pa_t^\al u(x,t)=\pa_x^2u(x,t), & 0<x<1,\ 0<t<T,\\
-\pa_x u(0,t)+\rho(t)u(0,t)=h_0(t), & 0<t<T,\\
\pa_x u(1,t)+\rho(t)u(1,t)=h_1(t), & 0<t<T,\\
u(x,0)=a(x), & 0<x<1,
\end{cases}
\]
the papers Wei and Wang \cite{WW16} and Wei and Zhang \cite{WZ16} discuss an inverse problem of determining $\rho(t)$, $0<t<T$ by data $u(0,t)$, $0<t<T$, and study numerical methods by reducing the inverse problem to a Volterra equation in $\rho$.

\section{Concluding remarks}\label{sec-con}

Taking the simple formulation \eqref{eq-ibvp-u} as a model problem, in this chapter we mainly reviewed the theoretical results on the determination of temporal and spatial components in source terms by several kinds of observation data. It reveals that the fractional derivative in time results in essential difficulties in treating these problems more than classical partial differential equations, so that the available arguments are limited, for example, representation formulae of the solution to the initial-boundary value problem \eqref{eq-ibvp-u} by the Mittag-Leffler functions. Thus we should recognize wider varieties of works on inverse problems for fractional partial differential equations.

In this chapter, our survey concentrates on theoretical works, but by natural necessity various works on numerical methods for related inverse problems have been continuously published. Numerical researches for inverse problems for fractional differential equations are tremendously expanding and here we are restricted to making a partial list of related works: Chi, Li and Jia \cite{CLJ11}, Jin and Rundell \cite{JR15}, Murio and Mej\'ia \cite{MM08}, Tian, Li, Deng and Wu \cite{TLDW12}, Wang and Wei \cite{WW15}, Wang, Yamamoto and Han \cite{WYH13}, Wei, Chen, Sun and Li \cite{WCSL10}, Wei and Wang \cite{WW14a}, Wei and Zhang \cite{WZ13}, Yang, Fu and Li \cite{YFL15}, Zhang, Li, Jia and Li \cite{ZLJL13}, Zhang and Wei \cite{ZW13}.

Finally, we mention several prospects on future topics. Compared with problems with order $\al\in(0,1)$, less has been done for the cases of $\al\in(1,2)$ which may also have practical significance. On the other hand, for realistic applications, other kinds of sources should also be taken into account, for example, the multiple point source $\sum_{i=1}^N\rho_i(t)\de(x-x_i)$ and the moving source $g(x-\xi(t))$. Numerically, it is preferable to develop advanced regularization methods capturing the fractional essence instead of direct optimization techniques.

\bibliographystyle{unsrt}

\begin{thebibliography}{1}

\bibitem{A75}
R. A. Adams, {\it Sobolev Spaces}, Academic Press, New York, 1975.

\bibitem{AKM13}
T. Aleroev, M. Kirane, and S. Malik, Determination of a source term for a time fractional diffusion equation with an integral type over-determining condition, {\it Electron. J. Differ. Equ.}, {\bf270} (2013), 1--16.

\bibitem{BY17}
M. Bellassoued and M. Yamamoto, {\it Carleman Estimates and Applications to Inverse Problems for Hyperbolic Systems}, Springer Japan, Tokyo, 2017.

\bibitem{CLN13}
J. Cheng, C.-L. Lin, and G. Nakamura, Unique continuation property for the anomalous diffusion and its application, {\it J. Differ. Equ.}, {\bf254} (2013), 3715--3728.

\bibitem{CLJ11}
G. Chi, G. Li, and X. Jia, Numerical inversions of a source term in the FADE with a Dirichlet boundary condition using final observations, {\it Comput. Math. Appl.}, {\bf62} (2011), 1619--1626.

\bibitem{CY96}
M. Choulli and M. Yamamoto, Generic well-posedness of an inverse parabolic problem--the H\"older-space approach, {\it Inverse Probl.}, {\bf12} (1996), 195--205.

\bibitem{CY97}
M. Choulli and M. Yamamoto, An inverse parabolic problem with non-zero initial condition, {\it Inverse Probl.}, {\bf13} (1997), 19--27.

\bibitem{CY99}
M. Choulli and M. Yamamoto, Generic well-posedness of a linear inverse parabolic problem with diffusion parameters, {\it J. Inverse Ill-Posed Probl.}, {\bf7} (1999), 241--254.

\bibitem{EK04}
S. D. Eidelman and A. N. Kochubei, Cauchy problem for fractional diffusion equations, {\it J. Differ. Equ.}, {\bf199} (2004), 211--255.

\bibitem{FK16}
K. Fujishiro and Y. Kian, Determination of time dependent factors of coefficients in fractional diffusion equations, {\it Math. Control Relat. Fields}, {\bf6} (2016), 251--269.

\bibitem{GLY15}
R. Gorenflo, Yu. Luchko, and M. Yamamoto, Time-fractional diffusion equation in the fractional Sobolev spaces, {\it Fract. Calc. Appl. Anal.}, {\bf18} (2015), 799--820.

\bibitem{H81}
D. Henry, {\it Geometric Theory of Semilinear Parabolic Equations}, Springer-Verlag, Berlin, 1981.

\bibitem{I91}
V. Isakov, Inverse parabolic problems with the final overdetermination, {\it Commun. Pure Appl. Math.}, {\bf44} (1991), 185--209.

\bibitem{JLLY17}
D. Jiang, Z. Li, Y. Liu, and M. Yamamoto, Weak unique continuation property and a related inverse source problem for time-fractional diffusion-advection equations, {\it Inverse Probl.}, {\bf33} (2017), 055013.

\bibitem{JR15}
B. Jin and W. Rundell, A tutorial on inverse problems for anomalous diffusion processes, {\it Inverse Probl.}, {\bf31} (2015), 035003.

\bibitem{K76}
T. Kato, {\it Perturbation Theory for Linear Operators}, Springer-Verlag, Berlin, 1976.

\bibitem{K18}
A. Kawamoto, H\"older stability estimate in an inverse source problem for a first and half order time fractional diffusion equation, {\it Inverse Probl. Imaging}, {\bf12} (2018), 315--330.

\bibitem{KM11}
M. Kirane and S. A. Malik, Determination of an unknown source term and the temperature distribution for the linear heat equation involving fractional derivative in time, {\it Appl. Math. Comput.}, {\bf218} (2011), 163--170.

\bibitem{KMA13}
M. Kirane, S. Malik, and M. Al-Gwaiz, An inverse source problem for a two dimensional time fractional diffusion equation with nonlocal boundary conditions, {\it Math. Methods Appl. Sci.}, {\bf36} (2013), 1056--1069.

\bibitem{KT04}
M. V. Klibanov and A. A. Timonov, {\it Carleman Estimates for Coefficient Inverse Problems and Numerical Applications},
VSP, Utrech, 2004.

\bibitem{KY17}
A. Kubica and M. Yamamoto, Initial-boundary value problems for fractional diffusion equations with time-dependent coeffcients, {\it Fract. Calc. Appl. Anal.}, {\bf21} (2018), 276--311.

\bibitem{LLY15}
Z. Li, Y. Liu, and M. Yamamoto, Initial-boundary value problems for multi-term time-fractional diffusion equations with positive constant coefficients, {\it Appl. Math. Comput.}, {\bf257} (2015), 381--397.

\bibitem{LN16}
C.-L. Lin and G. Nakamura, Unique continuation property for anomalous slow diffusion equation, {\it Commun. Partial Differ. Equ.}, {\bf41} (2016), 749--758.

\bibitem{LN2018}
C.-L. Lin and G. Nakamura, Unique continuation property for multi-terms time fractional diffusion equations, {\it Math. Ann.} (2 July 2018), https://doi.org/10.1007/s00208-018-1710-z.

\bibitem{LYY16}
J. J. Liu, M. Yamamoto, and L. L. Yan, On the reconstruction of unknown time-dependent boundary sources for time fractional diffusion process by distributing measurement, {\it Inverse Probl.}, {\bf32} (2016), 015009.

\bibitem{L17}
Y. Liu, Strong maximum principle for multi-term time-fractional diffusion equations and its application to an inverse source problem, {\it Comput. Math. Appl.}, {\bf73} (2017), 96--108.

\bibitem{LRY16}
Y. Liu, W. Rundell, and M. Yamamoto, Strong maximum principle for fractional diffusion equations and an application to an inverse source problem, {\it Fract. Calc. Appl. Anal.}, {\bf19} (2016), 888--906.

\bibitem{LZ17}
Y. Liu and Z. Zhang, Reconstruction of the temporal component in the source term of a (time-fractional) diffusion equation, {\it J. Phys. A}, {\bf50} (2017), 305203.

\bibitem{L09}
Yu. Luchko, Maximum principle for the generalized time-fractional diffusion equation, {\it J. Math. Anal. Appl.}, {\bf351} (2009), 218--223.

\bibitem{MM08}
D. A. Murio and C. E. Mej\'ia, Source terms identification for time fractional diffusion equation, {\it Rev. Colomb. Mat.}, {\bf42} (2008), 25--46.

\bibitem{P99}
I. Podlubny, {\it Fractional Differential Equations}, Academic Press, San Diego, 1999.

\bibitem{POV00}
A. I. Prilepko, D. G. Orlovsky, and I. A. Vasin, {\it Methods for Solving Inverse Problems in Mathematical Physics}, Marcel Dekkers, New York, 2000.

\bibitem{RW17}
Z. Ruan and Z. Wang, Identification of a time-dependent source term for a time fractional diffusion problem, {\it Appl. Anal.}, {\bf96} (2017), 1638--1655.

\bibitem{STY02}
S. Saitoh, V. K. Tuan, and M. Yamamoto, Reverse convolution inequalities and applications to inverse heat source problems,
{\it J. Inequal. Pure Appl. Math.}, {\bf3} (2002), 80.

\bibitem{SY11a}
K. Sakamoto and M. Yamamoto, Initial value/boundary value problems for fractional diffusionwave equations and applications to some inverse problems, {\it J. Math. Anal. Appl.}, {\bf382} (2011), 426--447.

\bibitem{SY11b}
K. Sakamoto and M. Yamamoto, Inverse source problem with a final overdetermination for a fractional diffusion equation, {\it Math. Control Relat. Fields}, {\bf1} (2011), 509--518.

\bibitem{TTU2016}
S. Tatar, R. Tinaztepe, and S. Ulusoy, Determination of an unknown source term in a space-time fractional diffusion equation, {\it J. Fract. Calc. Appl.}, {\bf 6} (2015), 83--90.

\bibitem{TU14}
S. Tatar and S. Ulusoy, S, An inverse source problem for a one-dimensional space-time fractional diffusion equation, {\it Appl. Anal.}, {\bf94} (2015), 2233--2244.

\bibitem{TLDW12}
W. Tian, C. Li, W. Deng, and Y. Wu, Regularization methods for unknown source in space fractional diffusion equation, {\it Math. Comput. Simul.}, {\bf85} (2012), 45--56.

\bibitem{TE94}
I. V. Tikhonov and Yu. S. Eidelman, Problems on correctness of ordinary and inverse problems for evolutionary equations of a special form, {\it Math. Notes}, {\bf56} (1994), 830--839.

\bibitem{T26}
E. C. Titchmarsh, The zeros of certain integral functions, {\it Proc. London Math. Soc.}, {\bf2} (1926), 283--302.

\bibitem{US07}
S. R. Umarov and E. M. Saidamatov, A generalization of Duhamels principle for differential equations of fractional order, {\it Dokl. Math.}, {\bf75} (2007), 94--96.

\bibitem{Wang-Wu14}
H. Wang and B. Wu, On the well-posedness of determination of two coefficients in a fractional integrodifferential equation, {\it Chin. Ann. Math.}, {\bf35B} (2014), 447--468.

\bibitem{WW15}
J.-G. Wang and T. Wei, Quasi-reversibility method to identify a space-dependent source for the time-fractional diffusion equation, {\it Appl. Math. Model.}, {\bf39} (2015), 6139--6149.

\bibitem{WZW13}
J.-G. Wang, Y.-B. Zhou, and T. Wei, Two regularization methods to identify a space-dependent source for the time-fractional diffusion equation, {\it Appl. Numer. Math.}, {\bf68} (2013), 39--57.

\bibitem{WYH13}
W. Wang, M. Yamamoto, and B. Han, Numerical method in reproducing kernel space for an inverse source problem for the fractional diffusion equation, {\it Inverse Probl.}, {\bf29} (2013), 095009.

\bibitem{WCSL10}
H. Wei, W. Chen, H. Sun, and X. Li, A coupled method for inverse source problem of spatial fractional anomalous diffusion equations, {\it Inverse Probl. Sci. Eng.}, {\bf18} (2010), 945--956.

\bibitem{WLL16}
T. Wei, X. L. Li, and Y. S. Li, An inverse time-dependent source problem for a time-fractional diffusion equation, {\it Inverse Problems}, {\bf32} (2016), 085003.

\bibitem{WSL16}
T. Wei, L. Sun, and Y. Li, Uniqueness for an inverse space-dependent source term in a multi-dimensional time-fractional diffusion equation, {\it Appl. Math. Lett.}, {\bf61} (2016), 108--113.

\bibitem{WW14a}
T. Wei and J. Wang, A modified quasi-boundary value method for an inverse source problem of the time-fractional diffusion equation, {\it Appl. Numer. Math.}, {\bf78} (2014), 95--111.

\bibitem{WW16}
T. Wei and J. Wang, Determination of Robin coefficient in a fractional diffusion problem, {\it Appl. Math. Model.}, {\bf40} (2016), 7948--7961.

\bibitem{WZ13}
T. Wei and Z. Q. Zhang, Reconstruction of a time-dependent source term in a time-fractional diffusion equation, {\it Eng. Anal. Bound. Elem.}, {\bf37} (2013), 23--31.

\bibitem{WZ16}
T. Wei and Z. Q. Zhang, Robin coefficient identification for a time-fractional diffusion equation, {\it Inverse Probl. Sci. Eng.}, {\bf24} (2016), 647--666.

\bibitem{WW14b}
B. Wu and S. Wu, Existence and uniqueness of an inverse source problem for a fractional integrodifferential equation, {\it Comput. Math. Appl.}, {\bf68} (2014), 1123--1136.

\bibitem{XCY11}
X. Xu, J. Cheng, and M. Yamamoto, Carleman estimate for a fractional diffusion equation with half order and application,
{\it Appl. Anal.}, {\bf90} (2011), 1355--1371.

\bibitem{Y09}
M. Yamamoto, Carleman estimates for parabolic equations and applications, {\it Inverse Probl.}, {\bf25} (2009), 123013.

\bibitem{YZ12}
M. Yamamoto and Y. Zhang, Conditional stability in determining a zeroth-order coefficient in a half-order fractional diffusion equation by a Carleman estimate, {\it Inverse Probl.}, {\bf28} (2012), 105010.

\bibitem{YFL15}
F. Yang, C. Fu, and X. Li, The inverse source problem for time-fractional diffusion equation: stability analysis and regularization, {\it Inverse Probl. Sci. Eng.}, {\bf23} (2015), 969--996.

\bibitem{Z09}
R. Zacher, Weak solutions of abstract evolutionary integro-differential equations in Hilbert spaces, {\it Funkc. Ekvacioj}, {\bf52} (2009), 1--18.

\bibitem{ZLJL13}
D. Zhang, G. Li, X. Jia, and H. Li, Simultaneous inversion for space-dependent diffusion coefficient and source magnitude in the time fractional diffusion equation, {\it J. Math. Res.}, {\bf5} (2013), 65--78.

\bibitem{ZX11}
Y. Zhang and X. Xu, Inverse source problem for a fractional diffusion equation, {\it Inverse Probl.}, {\bf27} (2011), 035010.

\bibitem{ZW13}
Z. Q. Zhang and T. Wei, Identifying an unknown source in time-fractional diffusion equation by a truncation method, {\it Appl. Math. Comput.}, {\bf219} (2013), 5972--5983.

\end{thebibliography}

\end{document}